\documentclass[12pt]{article}
\usepackage{mathrsfs}
\usepackage{epic,eepic,epsf,epsfig}
\usepackage{amsfonts,srcltx,mathrsfs}
\usepackage{amsmath}
\usepackage{amssymb}
\usepackage{amsbsy}
\usepackage{float}
\usepackage{graphicx}
\usepackage{amsfonts}
\usepackage{color}

\usepackage{url}
\newtheorem{lem}{Lemma}[section]
\newtheorem{theorem}[lem]{Theorem}

\newtheorem{cor}[lem]{Corollary}

\newtheorem{prob}[lem]{Problem}

\newtheorem{case}{Case}
\newtheorem{prop}[lem]{Proposition}

\def\r{\rho} \def\s{\sigma}   
  
 \def\Ga{\Gamma}

 \def\lg{\langle} \def\rg{\rangle}
\def\nd{\mathrel{\bigm|\kern-.7em/}}

\def\f{\noindent}

\def\Aut{\hbox{\rm Aut}}

\def\Aut{\hbox{\rm Aut}}

\def\demo{\f {\bf Proof.}\hskip10pt}

\newcommand{\qed}{\mbox{\raisebox{0.7ex}{\fbox{}}} \vspace{4truemm}}

\def\P{\mathcal {P}}

\begin{document}

\title{Four infinite families of chiral $3$-polytopes of type $\{4, 8\}$ with solvable automorphism groups}

\author{ \\ Dong-Dong Hou, Tian-Tian Zheng, Rui-Rui Guo\\
{\small Department of Mathematics, Shanxi Normal University}\\
{\small Taiyuan, Shanxi
041004, P.R. China}\\
}

\date{}
\maketitle

\footnotetext{E-mails: holderhandsome$@$bjtu.edu.cn, 709127994@qq.com, 13368420984@163.com
}
\begin{abstract}

We construct four infinite families of chiral $3$-polytopes of type $\{4, 8\}$,
with $1024m^4$, $2048m^4$, $4096m^4$ and $8192m^4$ automorphisms for every positive integer $m$, respectively.
The automorphism groups of these polytopes are solvable groups, and when $m$ is a power of $2$,
they provide examples with automorphism groups of order $2^n$ where $n \geq 10$. 
(On the other hand, no chiral polytopes of type $\{4, 8\}$ exist for $n \leq 9$.)
In particular,  our families give a partial answer to a problem proposed by Schulte and Weiss
in [Problems on polytopes, their groups, and realizations, {\em Period. Math. Hungar.} 53 (2006), 231-255]
and a problem proposed by Pellicer in [Developments and open problems on chiral polytopes, {\em Ars Math. Contemp} 5 (2012), 333-354].

\bigskip
\f {\bf Keywords:} Chiral $3$-polytope, Automorphism group, $2$-group.\\
\end{abstract}

%%%%%%%%%%%
\section{Introduction}
%%%%%%%%%%%

~\quad In~\cite{Problem}, Schulte and Weiss proposed the following problem:
\begin{prob}\label{problem1}
Characterize the groups of orders $2^n$, with $n$ a positive
integer, which are automorphism groups of regular or chiral
polytopes.
\end{prob}

Let $\P$ be a regular or chiral $d$-polytope, whose automorphism group has order $2^n$ with type $\{2^{k_1}, 2^{k_2}, \cdots, 2^{k_d}\}$. Here $k_i >1$ for $i \in \{1, 2, \cdots, d\}$.
The atlas~\cite{atles1} contains information about all regular or chiral polytopes with automorphism group of order at most 2000.
The first author, Feng and Lemmans in~\cite{HFL, HFL1} shows that  if $k_1+k_2+\cdots +k_d \leq n-1$, then there exists an
regular $d$-polytope of order $2^n$ with type  $\{2^{k_1}, 2^{k_2}, \cdots, 2^{k_d}\}$ for $n \geq 10$. That means all possible 
type can be achieved for regular polytope of order $2^n$. However, there are
just a few results for chiral polytopes of order $2^n$, see~\cite{CFH, Zhang'sPhdThesis}. 
Even for the case $d=3$, to the best of our knowledge, the only infinite family is the so-called tight chiral-polytopes, that is, it has type
$\{2^{k_1}, 2^{k_2}\}$ and has order $2^n=2^{k_1k_2}$, see~\cite{TightChiralPolyhedra}. On the other hand, by the famous book~\cite[Chapter 8]{BooksCoxeter}, one can see that $\{2^{k_1}, 2^{k_2}\} \ne \{4, 4\}$.
Inspired by these results listed above, we  naturally  consider the following problem:
\begin{prob}\label{problem2}
what is the smallest type of chiral $3$-polytopes of order $2^n$? 
\end{prob}

Here we construct four families of chiral $3$-polytopes with type $\{4, 8\}$.
Each family contains one example with $1024m^4$, $2048m^4$, $4096m^4$ or $8192m^4$ automorphisms, respectively, for every integer $m \geq 1$.
In particular, if we let $m$ be an arbitrary power of $2$, say $2^k$ (with $k \geq 0$), the the automorphism group has order $2^{10+4k}$, $2^{11+4k}$, $2^{12+4k}$ or $2^{13+4k}$,
which can be expressed as $2^n$ for an arbitrary integer $n \geq 10$. 
It means that the smallest type of chiral $3$-polytopes of order $2^n$ is $\{4, 8\}$ for $n \geq 10$.
Hence, these results givie an  answer to problem~\ref{problem2} and a partial answer to problem~\ref{problem1}.

On the other hand, chiral $3$-polytopes are also known as chiral maps.  The genus of a chiral map is the genus of the carrier surface. 
The genus $g$ and type  $\{k_1, k_2\}$ of a chiral map are related via the Euler–Poincar\'e formula by
$2-2g = \chi = |G|(\frac{1}{k_1}+\frac{1}{k_2}-\frac{1}{2})$, where $\chi$ is the Euler characteristic of the carrier surface
and $G$ is the automorphism group of chiral maps. All chiral maps on the tours was
described by Coxeter in~\cite{BooksCoxeter}.
In fact, no chiral maps lies on a non-orientable surface and there are no chiral 
maps on orientable surfaces from genus $2$ to genus $6$. 
 These resultes were extended up to genus $15$ by Conder and Dobcs\'anyi, with
the help of computational methods, and then more recently 
by Conder much further, up to genus $301$, see~\cite{Conder+EulerCharacteristic,Conder+date}.
In~\cite{ProblemOfChiralPolytopes}, Pellicer proposed the following problem:
\begin{prob}\label{problem3}
Determine all positive integers $g$ for which there are chiral maps(polyhedra) on orientable surfaces with genus $g$.
\end{prob}

It was proved in~\cite{Conder+Siran+Tucker} that orientable surfaces with genera $p+1$, ($\chi=-2p$), with $p$ prime and $p-1$ not divisible by $3, 5$ or $8$,
do not admit chiral maps. Our chiral polytopes shows that for every $n \geq 10$, there exists a chiral polytope(map) of type $\{4, 8\}$ on
an orientable surfaces with genera $g=2^{n-4}+1\geq  65$, ($\chi=-2^{n-3}$), with $2^n$ automorphisms. 
Hence, these results given a partial answer to problem~\ref{problem3}.

Our main result is the following theorem:

\begin{theorem}\label{maintheorem}
For every  positive integer $m \geq 1$, there exist chiral $3$-polytopes $\P^{1}_m, \P^{2}_m$, $\P^{3}$ and $\P^{4}_m$ of type $\{4,8\}$
with solvable automorphism groups of order $1024m^4, 2048m^4$, $4096m^4$ and $8192m^4$, respectively.
\end{theorem}

As a special case we have the following Corollary, which is an immediate consequence of Theorem~\ref{maintheorem}
when $m$ is taken as a power of $2$.

\begin{cor}
For $n \geq 10$, the smallest type of chiral $3$-polytopes of order $2^n$ is $\{4, 8\}$, and the corresponding surfaces has 
Euler characteristic $\chi=-2^{n-3}$ $(g=2^{n-4}+1)$.
\end{cor}

Together with Conder's Date\cite{atles1,Conder+date}, we have the following Corollary.

\begin{cor}
For any $l \geq 0$, let $\mathcal{S}$ be an orientable surface with Euler characteristic $\chi=-2^{l}$, then there exists a chiral map on $\mathcal{S}$
if and only if $l \ne 1, 2, 3$ and $4$, that is $\chi \ne -2, -4, -8$ and $-16$.
\end{cor}

%%%%%%%%%%%%%%%
\section{Additional background}
\label{background}
%%%%%%%%%%%%%%%

In this section we give some further background that may be helpful for the rest of the paper, see~\cite{ARP,ChiralPolytopes,GroupBook}.

\subsection{Abstract polytopes: definition, structure and properties}

An abstract polytope of rank $n$ is a partially ordered set $\P$ endowed with a strictly monotone rank function
with range $\{-1, 0, \cdots, n\}$, which satisfies four conditions, to be given shortly.

The elements of $\P$ are called \emph{faces} of $\P$. More specifically, the elements of $\P$ of rank $j$ are called $j$-faces,
and a typical $j$-face is denoted by $F_j$.
Two faces $F$ and $G$ of $\P$ are said to be \emph{incident} with each other if $F \leq G$ or $F \geq G$ in $\P$.
A \emph{chain} of $\P$ is a totally ordered subset of $\P$, and is said to have \emph{length} $i$ if it contains exactly $i+1$ faces.
The maximal chains in $\P$ are called the \emph{flags} of $\P$.
Two flags are said to be $j$-\emph{adjacent} if they differ in just one face of rank $j$,
or simply \emph{adjacent} (to each other) if they are $j$-adjacent for some~$j$.
% in which case they may be denoted by $\Phi$ and $\Phi^j$.
Also if $F$ and $G$ are faces of $\P$ with $F \leq G$, then the set $\{\,H \in \P \ |\ F \leq H \leq G\,\}$ is
called a \emph{section} of $\P$, and is denoted by $G/F$.  Such a section has rank $m-k-1$,
where $m$ and $k$ are the ranks of $G$ and $F$ respectively.
A section of rank $d$ is called a $d$-section.

\smallskip
We can now give the four conditions that are required of $\P$ to make it an abstract polytope.
These are listed as (P1) to (P4) below:
\begin{itemize}
\item [(P1)]  $\P$ contains a least face and a greatest face, denoted by $F_{-1}$ and $F_n$, respectively.
\item [(P2)]  Each flag of $\P$ has length $n+1$ (so has exactly $n+2$ faces, including $F_{-1}$ and $F_n$).
\item [(P3)]  $\P$ is \emph{strong flag-connected}, which means that any two flags $\Phi$ and $\Psi$ of $P$ can
be joined by a sequence of successively adjacent flags $\Phi=\Phi_{0}, \Phi_1, \cdots, \Phi_k=\Psi$, each of which contains $\Phi \cap \Psi$.
\item [(P4)]  The rank $1$ sections of $\P$ have a certain homogeneity property known as the \emph{diamond condition},
namely as follows: if $F$ and $G$ are incidence faces of $\P$, of ranks $i-1$ and $i+1$, respectively, where $0 \le i \le n-1$,
then there exist precisely \emph{two} $i$-faces $H$ in $\P$ such that $F< H< G$.
\end{itemize}

\noindent
An easy case of the diamond condition occurs for polytopes of rank 3 (or polyhedra): if $v$ is a vertex of same face $f$,
then there are two edges that are incident with both $v$ and $f$.

\smallskip
If $F_{n-1}$ is a facet (of rank $n-1$), then
the section $F_{n-1}/F_{-1}$ is also called a {\em facet} of $\P$,
while if $F_0$ is a vertex, then the section
$F_{n}/F_{0} = \{G \in \P \ |\ F_{0} \leq G\}$ is called a {\em vertex-figure} of $\P$ at $F_0$.
Next, every $2$-section $G/F$ of $\P$ is isomorphic to the face lattice of a polygon.
Now if it happens that the number of sides of every such polygon depends only on the rank of $G$,
and not on $F$ or $G$ itself, then we say that the polytope $\P$ is {\em equivelar}.
In this case, if $k_i$ is the number of edges of every $2$-section between an $(i-2)$-face and an $(i+1)$-face of $\P$,
for $1 \leq i \leq n$, then the expression $\{k_1, k_2, \cdots, k_{n-1}\}$ is called the Schl\"afli type of $\P$.
(For example, if $\P$ has rank 3, then $k_1$ and $k_2$ are the valency of each vertex and the size of each face,
respectively.)

\subsection{Automorphisms of polytopes}

An \emph{automorphism} of an abstract polytope $\P$ is an order-preserving permutation of its elements.
In particular, every automorphism preserves the set of faces of any given rank.
Under permutation composition, the set of all automorphisms of $\P$ forms a group, called the automorphism
group of $\P$, and denoted by $\Aut(\P)$ or sometimes more simply as $\Ga(\P)$.
Also it is not difficult to use the diamond condition and strong flag-connectedness to prove that if an automorphism
preserves of flag of $\P$, then it fixes every flag of $\P$ and hence every element of $\P$.
It follows that $\Ga(\P)$ acts semi-regularly on flags of $\P$.

\smallskip
A polytope $\P$ is said to be \emph{regular} if its automorphism group $\Ga(\P)$ acts transitively
(and hence regularly) on the set of flags of $\P$.
In this case, the number of automorphisms of $\P$ is as large as possible, and equal to the number of flags of $\P$.
In particular, $\P$ is equivelar, and the stabiliser in $\Ga(\P)$ of every 2-section of $\P$ induces the full dihedral group
on the corresponding polygon.
Moreover, for a given flag $\Phi$ and for every $i \in \{0,1,\dots,n-1\}$, the polytope $\P$ has a unique
automorphism $\rho_i$ that takes $\Phi$ to the unique flag $(\Phi)^{i}$ that differs from $\Phi$ in precisely its $i$-face,
and then the automorphisms $\rho_0,\rho_1,\dots,\rho_{n-1}$ generate $\Ga(\P)$ and satisfy the defining
relations for the string Coxeter group $[k_1, k_2, \cdots, k_{n-1}]$, where the $k_i$ are as given in the
previous subsection for the Schl\"afli type of $\P$.
Here, the {\em string Coxeter group} $[k_1, k_2, \cdots, k_{n-1}]$, is defined as the group with presentation

$$\lg \r_0, \r_1, \cdots, \r_{n-1} \ |\ \r_i^2=1 \ {\rm for} \ 0 \leq i \leq n-1, (\r_i\r_{i+1})^{k_{i+1}}=1 \ {\rm for} \ 0 \leq i \leq n-2,$$
$$(\r_i\r_j)^2=1 \ {\rm for} \ 0 \leq i < j-1  < n-1\rg.$$

They also satisfy a certain `intersection condition', which follows from the diamond and strong flag-connectedness
conditions.
These and many more properties of regular polytopes may be found in~\cite{ARP}.

\smallskip
We now turn to chiral polytopes, for which two good references are~\cite{ChiralPolytopes}  and~\cite{ProblemOfChiralPolytopes}.

\smallskip
A polytope $\P$ said to be {\em chiral\/} if its automorphism group $\Ga(\P)$ has two orbits on flags,
with every two adjacent flags lying in different orbits.
(Another way of viewing this definition is to consider $\P$ as admitting no `reflecting' automorphism
that interchanges a flag with an adjacent flag.)
Here the number of flags of $\P$ is $2|\Ga(\P)|$,  and $\Ga(\P)$ acts regularly on each of two orbits.
Again $\P$ is equivelar, with the stabiliser in $\Ga(\P)$ of every 2-section of $\P$ inducing the full cyclic group on the corresponding polygon.

 For a given flag $\Phi$, denote by $F_i$ be the $i$-face in $\Phi$ for each $0\leq i\leq n$, and for every $j \in \{1,2,\dots,n-1\}$, the chiral polytope $\P$ admits an automorphism $\s_j$ that fixes each $F_i$ with $i\not=j-1,j$ and cyclically permutes
consecutive $j$- and $(j-1)$-faces in the $2$-section $F_{j+1}/F_{j-2}$, 
that is, $\s_j$ takes $\Phi$ to the flag $(\Phi)^{j, j-1}$ which differs from $\Phi$ in precisely its $(j-1)$-and $j$-faces.
This automorphism $\sigma_j$ is the analogue of the abstract rotation $\r_{j-1}\r_{j}$ in the
regular case, for each $j$. These automorphisms $\s_1,\s_2,\dots,\s_{n-1}$ generate $\Ga(\P)$, and if $\P$ has Schl\"afli type $\{k_1, k_2, \dots, k_{n-1}\}$,
then they satisfy the defining relations for the orientation-preserving subgroup of (index $2$ in) the
string Coxeter group $[k_1, k_2, \cdots, k_{n-1}]$.
Also they satisfy a `chiral' form of the intersection condition, which is a variant of the one mentioned earlier for regular polytopes.

Chiral polytopes occur in pairs (or {\em enantiomorphic\/} forms), such that each member of the pair is the `mirror image' of the other.
Suppose one of them is $\P$, and has Schl\"afli type $\{k_1, k_2, \cdots, k_{n-1}\}$. Then $\Ga(\P)$ is isomorphic
to the quotient of the orientation-preserving subgroup $\Lambda^{\rm o}$ of the string Coxeter group $\Lambda = [k_1, k_2, \cdots, k_{n-1}]$
via some normal subgroup $K$.  By chirality, $K$ is not normal in the full Coxeter group $\Lambda$, but is conjugated by
any orientation-reversing element $c \in \Lambda$ to another normal subgroup $K^c$ which is the kernel of an epimorphism from $\Lambda^{\rm o}$
to the automorphism group $\Ga(\P^c)$ of the mirror image $\P^c$ of $\P$.

The automorphism groups of $\P$ and $\P^c$ are isomorphic to each other, but their canonical generating sets
satisfy different defining relations.   In fact, replacing the elements $\s_1$ and $\s_2$ in the canonical generating
tuple $(\s_1,\s_2,\s_3,\dots,\s_{n-1})$ by $\s_1^{-1}$ and $\s_1^{\,2}\s_2$ gives a set of generators for $\Ga(\P)$
that satisfy the same defining relations as a canonical generating tuple for $\Ga(\P^c)$,
but chirality ensures that there is no automorphism of $\Ga(\P)$ that takes $(\s_1,\s_2)$ to $(\s_1^{-1},\s_1^{\,2}\s_2)$
and fixes all the other $\s_j$.

Conversely, any finite group $G$ that is generated by $n-1$ elements $\s_1,\s_2,\dots,\s_{n-1}$ which satisfy both
the defining relations for $\Lambda^{\rm o}$ and the chiral form of the intersection condition is the `rotation subgroup' of an
abstract $n$-polytope $\P$ that is either regular or chiral.
Indeed, $\P$ is regular if and only if $G$ admits a group automorphism $\r$ of order $2$
that takes $(\s_1,\s_2,\s_3,\dots,\s_{n-1})$ to $(\s_1^{-1},\s_1^{\,2}\s_2,\s_3,\dots,\s_{n-1})$.

\smallskip
We now focus our attention on the rank $3$ case.
Here the generators $\s_1, \s_2$  for $\Ga(\P)$ satisfy the canonical relations
$\,\s_1^{k_1} = \s_2^{k_2} = (\s_1\s_2)^2  = 1$,
and the chiral form of the intersection condition can be abbreviated to
$\lg \s_1\rg \cap \lg \s_2 \rg =\{1\}$.

\subsection{Group theory}

We use  standard notation for group theory, as in~\cite{GroupBook} for example.
We also need the following, which are elementary and so we give them without proof.

\begin{prop}\label{freeabelian}
Let $G$ be the free abelian group $\mathbb{Z}\oplus \mathbb{Z}\oplus \mathbb{Z}\oplus \mathbb{Z}$ of rank $4$, generated by four elements $x_1, x_2, x_3$ and $x_4$
subject to the single defining relation $[x_i,x_j] = 1, i, j \in \{1,2,3,4\}$.
Then for every positive integer $m$, the subgroup $G_m=\lg x_1^m, x_2^m, x_3^m, x_4^m\rg$ is characteristic in $G$, with index $|G:G_m|=m^4$.
\end{prop}

Finally, we will use some Reidemeister-Schreier theory, which produces a defining presentation
for a subgroup $H$ of finite index in a finitely-presented group $G$.
An easily readable reference for this is~\cite[Chapter IV]{Johnson}, but in practice we use its implementation
as the {\tt Rewrite} command in the {\sc Magma} computation system~\cite{BCP97}.
We also found the groups that we use in the next section with the help of {\sc Magma} in constructing
and analysing some small examples.

\section{Proof of Theorem~\ref{maintheorem}}\label{Maim results}

\demo
We begin by defining ${\cal U}$ as the finitely-presented group
$$
\lg\, a, b \ | \ a^4 = b^8 = (ab)^2 =[a^2,b^2]^2=(ab^3a^2b^4)^2= 1 \,\rg.
$$
This group ${\cal U}$ has four normal subgroups of index $1024, 2048, 4096$ and $8192$.
The quotients of ${\cal U}$ by each of these give the initial members of our four infinite families.

{\bf Case 1:} Take $N^1$ as the subgroup of ${\cal U}$ generated by $x_1, x_2, x_3, x_4$, where
$$ x_1=(b^{-2}a)^4,  x_2=(b^4)^{ab^{-1}}(b^4)^{a^{-1}},$$
$$ x_3=(b^2a^2)^4,~~~ x_4=((b^2a^2)^4)^a.$$
A short computation with {\sc Magma} shows that $N^1$ is normal in ${\cal U}$, with index $1024$. In fact,
the defining relations for ${\cal U}$ can be used to show that
\begin{center}
\begin{tabular}{llll}
 $x_1^a = x_2x_4^{-1}$ \quad &  $x_2^a = x_1^{-1}x_3^{-1}$ \quad &  $x_3^a = x_4$ \quad  &    and \quad $x_4^a = x_1^{-1}x_2^{-1}$,  \\[+2pt]
 $x_1^b = x_1x_4$ \quad &  $x_2^b = x_1^{-1}$ \quad&  $x_3^b = x_1^{-1}x_2x_3^{-1}x_4^{-1}$ \quad& and \quad $x_4^b = x_3$.  \\[-3pt]
\end{tabular}
\end{center}
One way to prove these relations is by hand, which we leave as a challenging exercise for the interested reader.
Another is by  a partial enumeration of cosets of the identity subgroup in ${\cal U}$.
For example, if this is done using the {\tt ToddCoxeter} command in {\sc Magma}, allowing the definition of just 110000 cosets,
then multiplication by each of the words $x_1^a(x_2x_4^{-1})^{-1}$ and $x_1^b(x_1x_4)^{-1}$ is found to fix the trivial coset,
and therefore $x_1^a(x_2x_4^{-1})^{-1}= 1 =x_1^b(x_1x_4)^{-1}$.

Also {\sc Magma}'s {\tt Rewrite} command gives a defining presentation for $N^1$, with 
$$[x_1,x_2] = [x_1, x_3] = [x_1, x_4]=[x_2, x_3]=[x_2,x_4]=[x_3, x_4]=1.$$
Hence the normal subgroup $N^1$ is free abelian of rank 4.

\smallskip
The quotient ${\cal U}/N$ is isomorphic to the automorphism group of the chiral 3-polytope
of type $\{4,8\}$ with 1024 automorphisms listed at~\cite{atles1}.

\smallskip
Now for any positive integer $m$, let $N^1_m$ be the subgroup generated by $x_1^m, x_2^m, x_3^m$
and $x_4^m$. By Proposition~\ref{freeabelian}, we know that $N_m$ is characteristic in $N$
and hence normal in ${\cal U}$, with index $|\,{\cal U}:N_m| = |\,{\cal U}:N||N:N_m| = 1024m^4$.
Moreover, in the quotient $G^1_m = {\cal U}/N_m$, the subgroup  $N/N_m$ is abelian and normal,
with quotient $({\cal U}/N_m)/(N/N_m) \cong  {\cal U}/N$ being a $2$-group, and so $G^1_m$ is solvable.

Next, we claim that $\lg a\rg \cap \lg b \rg =\{1\}$ in $G^1_m$. Since $a^4=1$, we have $|\lg a\rg \cap \lg b \rg| =1, 2$ or $4$.
If $|\lg a\rg \cap \lg b \rg|=4$, then $\lg a\rg \leq \lg b\rg$ and $G_m=\lg b\rg$,  contradicting $|G^1_m|=1024m^4$.
If $|\lg a\rg \cap \lg b\rg| =2$, then $\lg a^2 \rg =\lg b^4\rg$. It follows that $\lg a^2\rg \unlhd G^1_m$. Consider the quotient
group $G^1_m/\lg a^2\rg$. It is easy to see $(a\lg a^2\rg)^2=(b\lg a^2\rg)^4=(ab\lg a^2\rg)^2=\lg a^2\rg$. Then $|G^1_m/\lg a^2\rg| \leq 8$
and hence $|G^1_m| = |G^1_m/\lg a^2\rg|\cdot |\lg a^2\rg| \leq 8\cdot 2=16$, which is impossible because $|G^1_m|=1024m^4$.

Finally, if there exists an automorphism $\rho$ of $G^1_m$
taking $(a, b)$ to $(a^{-1}, a^2b)$, and  it  follows from the
relation $1 = (ab^3a^2b^4)^2$ that also
$1 =(a^{-2}(a^2b)^{3}a^{-2}(a^2b)^{4})^2$ in $G^1_m$. 
However, by {\sc Magma}, $a^{-2}(a^2b)^{3}a^{-2}(a^2b)^{4}$ has order $4$ in $G^1_1$ $(m=1)$.
On the other hand, $G^1_1$ is the quotient group of $G^1_m$. Then $a^{-1}b^{-3}a^{-2}b^{-4}$ has order at least $4$ in $G^1_m$,
which is impossible because $a^{-2}(a^2b)^{3}a^{-2}(a^2b)^{4}$ has order $2$ in $G^1_m$.

\smallskip
Thus $\P^1_m$ is chiral, with solvable automorphism group $G^1_m$ of order $1024m^4$.

\f {\bf {\sc \bf Magma} programs for Case 1:} Here we provide some {\sc Magma} code for determining
$\mathcal{U} \unrhd N^1 = \mathbb{Z}  \oplus \mathbb{Z}  \oplus\mathbb{Z}  \oplus\mathbb{Z} $,
$x_1^a = x_2x_4^{-1}$, and $x_1^b = x_1x_4$.

{\tt
\medskip

U<a,b>:=Group<a,b|a$^4$,b$^8$,(a*b)$^2$,(a$^2$,b$^2$)$^2$,(a*b$^3$*a$^2$*b$^4$)$^2$>;
\medskip

x1:=(b$^{-2}$*a)$^4$;
\medskip

 x2:=(b$^4)^{{\tt a*b^{-1}}}$*(b$^4)^{{\tt a^{-1}}}$;
\medskip

 x3:=(b$^2$*a$^2$)$^4$;
\medskip

 x4:=((b$^2$*a$^2$)$^4$)$^{\tt a}$;

\medskip

 N$^{1}$:=sub<U | x1, x2, x3, x4>;

\medskip

 Index(U, N$^{1}$); // {\rm Returns the index of N$^{1}$ in U}

\medskip

 IsNormal(U,N$^{1}$);  // {\rm Return true if N$^{1}$ is a normal subgroup of U and false otherwise}

\medskip

AbelianQuotientInvariants(N$^{1}$); // {\rm This function computes the elementary divisors of the derived quotient group N$^{1}$/[N$^{1}$, N$^{1}$]}

\medskip

 Rewrite(U, N$^{1}$); // {\rm Compute a defining set of relations for N$^{1}$ on the existing generators}

\medskip

 nc,ctb:=ToddCoxeter(U,sub<U|>:CosetLimit:=110000); // {\rm  This function attempts to build up a coset table of $1$ in U using the Todd-Coxeter procedure}

\medskip

 for g in [a,b] do 
\medskip

~~~~for i in [-1,0,1] do 
\medskip

~~~~~~~~for j in [-1,0,1] do 
\medskip

~~~~~~~~~~~~for k in [-1,0,1] do 
\medskip

~~~~~~~~~~~~~~~~for l in [-1,0,1] do 
\medskip

 ~~~~~~~~~~~~~~~~~~~~rel:=(x1$^{\tt g}$)*(x1$^{\tt i}$*x2$^{\tt j}$*x$3^{\tt k}$*x4$^{\tt l}$)$^{{\tt -1}}$;

\medskip

~~~~~~~~~~~~~~~~~~~~~~~~if ctb(1,rel) eq 1 then
\medskip

~~~~~~~~~~~~~~~~~~~~~~~~print g,":",i,j,k,l;

\medskip

 ~~~~~~~~~~~~~~~~~~~~end if;
\medskip

 ~~~~~~~~~~~~~~~~end for;
\medskip

~~~~~~~~~~~~end for;
\medskip

 ~~~~~~~~end for;
\medskip

 ~~~~end for;

 end for;}

\begin{figure}
 
 \center{\includegraphics[width=20cm]  {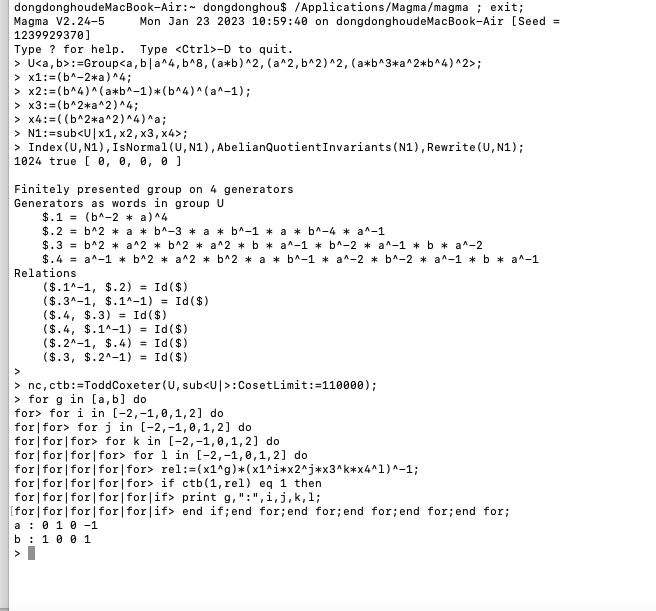}} 
 
 \caption{{\sc \bf Magma} programs for Case 1} 
 
 \end{figure}

\smallskip
{\bf Case 2:} Take $N^2$ as the subgroup of ${\cal U}$ generated by $y_1, y_2, y_3, y_4$, where
$$ y_1=b^2a^2ba^{-1}b^4a^{-1}b^2a^{-1}ba^{-1}, y_2=(b^{-1}a)^8, $$
$$y_3=b^3a^2b^2ab^{-1}ab^{-1}ab^{-3}a^{-1}, y_4=(ab^{-3})^4.$$
A short computation with {\sc Magma} shows that $N^2$ is normal in ${\cal U}$, with index $2048$,
and moreover, the {\tt Rewrite} command tells us that $N^2$ is free abelian of rank $4$. In this case,
the defining relations for ${\cal U}$ give
\begin{center}
\begin{tabular}{llll}
 $y_1^a = y_3$ \quad &  $y_2^a = y_1y_3^{-1}y_4^{-1}$ \quad &  $y_3^a = y_3y_4$ \quad  &    and \quad $y_4^a = y_2y_4^{-1}$,  \\[+2pt]
 $y_1^b = y_2^{-1}y_3^{-1}$ \quad &  $y_2^b = y_1^{-1}y_3^{-1}y_4^{-1}$ \quad&  $y_3^b = y_1$ \quad& and \quad $y_4^b = y_1^{-1}y_3^{-1}$.  \\[-3pt]
\end{tabular}
\end{center}

The quotient ${\cal U}/N^2$ is isomorphic to the automorphism group of the chiral 3-polytope
of type $\{4,8\}$ with 2048 automorphisms listed at~\cite{atles1}.

Now for any positive integer $m$, let $N^2_m$ be the subgroup generated by $y_1^m, y_2^m, y_3^m$, and $y_4^m$. Using
Proposition~\ref{freeabelian}, we find that $N^2_m$  is characteristic in $N_2$ and hence normal in ${\cal U}$, with 
index $|{\cal U}:N^2_m|=|{\cal U}: N^2| |N^2: N^2_m|=2048m^4$. Also the quotient $G^2_m={\cal U}/N^2_m$ is solvable.

Moreover, the same argument as used in case(1) shows that $\mathcal{P}^2_m$ is chiral. Thus, $\mathcal{P}^2_m$
is chiral, with automorphism group $G^2_m$ of order $2048m^4$.

{\bf Case 3:} Take $N^3$ as the subgroup of ${\cal U}$ generated by $z_1, z_2, z_3, z_4$, where
$$z_1=(b^{-1}a)^8, z_2=(b^{-3}a)^{4},$$
$$ z_3=(ab^{-1})^8, z_4=(ab^2a^{-1}ba^{-1})^4.$$
A short computation with {\sc Magma} shows that $N^3$ is normal in ${\cal U}$, with index $4096$,
and moreover, the {\tt Rewrite} command tells us that $N^3$ is free abelian of rank $4$. In this case,
the defining relations for ${\cal U}$ give
\begin{center}
\begin{tabular}{llll}
 $z_1^a =z_3^{-1}$ \quad &  $z_2^a =z_2^{-1}z_3^{-1} $ \quad &  $z_3^a = z_1$ \quad  &    and \quad $z_4^a =z_1^{-1}z_4$,  \\[+2pt]
 $z_1^b = z_2z_4^{-1}$ \quad &  $z_2^b =z_1^{-1}z_2 $ \quad&  $z_3^b =z_1$ \quad& and \quad $z_4^b =z_3^{-1}z_4^{-1}$.  \\[-3pt]
\end{tabular}
\end{center}

Now for any positive integer $m$, let $N^3_m$ be the subgroup generated by $z_1^m, z_2^m, z_3^m$, and $z_4^m$. Using
Proposition~\ref{freeabelian}, we find that $N^3_m$  is characteristic in $N^3$ and hence normal in ${\cal U}$, with 
index $|{\cal U}:N^3_m|=|{\cal U}: N^3| |N^3: N^3_m|=4096m^4$. Also the quotient $G^3_m={\cal U}/N^3_m$ is solvable.

Moreover, the same argument as used in case(1) shows that $\mathcal{P}^3_m$ is chiral. Thus, $\mathcal{P}^3_m$
is chiral, with automorphism group $G^3_m$ of order $4096m^4$.

{\bf Case 4:} Take $N^4$ as the subgroup of ${\cal U}$ generated by $w_1, w_2, w_3, w_4$, where
$$w_1=(ab^{-1})^8, w_2=((b^{-1}a)^8)^b,w_3=((b^{-2}a)^8)^{b},  $$
$$w_4=((a^2b^2)^4)^{b^{-1}}((b^{-2}a^2)^4)^a.$$
A short computation with {\sc Magma} shows that $N^4$ is normal in ${\cal U}$, with index $8192$,
and moreover, the {\tt Rewrite} command tells us that $N^4$ is free abelian of rank $4$. In this case,
the defining relations for ${\cal U}$ give
\begin{center}
\begin{tabular}{llll}
 $w_1^a =w_2 w_3^{-1}w_4^{-1}$ \quad &  $w_2^a = w_2w_3^{-1}$ \quad &  $w_3^a = w_3^{-1}$ \quad  &    and \quad $w_4^a = w_1w_2$,  \\[+2pt]
 $w_1^b = w_2 w_3^{-1}w_4^{-1}$ \quad &  $w_2^b = w_2^{-1}w_3$ \quad&  $w_3^b = w_1^{-1}w_2^{-1}w_3$ \quad& and \quad $w_4^b = w_1w_2^{-1}$.  \\[-3pt]
\end{tabular}
\end{center}

The quotient ${\cal U}/N^4$ is isomorphic to the automorphism group of the chiral 3-polytope
of type $\{4,8\}$ with 8192 automorphisms.

Now for any positive integer $m$, let $N^4_m$ be the subgroup generated by $w_1^m, w_2^m, w_3^m$, and $w_4^m$. Using
Proposition~\ref{freeabelian}, we find that $N^4_m$  is characteristic in $N_4$ and hence normal in ${\cal U}$, with 
index $|{\cal U}:N^4_m|=|{\cal U}: N^4| |N^4: N^4_m|=8192m^4$. Also the quotient $G^4_m={\cal U}/N^4_m$ is solvable.

Moreover, the same argument as used in case(1) shows that $\mathcal{P}^4_m$ is chiral
(in this case, $a^{-2}(a^2b)^{3}a^{-2}(a^2b)^{4}$ has order $8$ in $G^1_4$ $(m=1)$.). Thus, $\mathcal{P}^4_m$
is chiral, with automorphism group $G^4_m$ of order $8192m^4$.  \hfill\qed

\medskip
\f {\bf Acknowledgements:}
The authors acknowledge 
significant use of the Magma system~\cite{BCP97} in helping
us to make and test many of the observations described in this paper and the support from  National Natural Science Foundation of China (12201371),
Fundamental Research Program of Shanxi Province (20210302124078), Scientific and Technologial Innovation Programs of
Hight Education Institutions in Shanxi (2021L254).
\\[-12pt]

%%%%%%%%%%%%%%

\end{document}